\documentclass[twoside,a4paper,11pt,centertags]{amsart}
\usepackage{amsmath,amssymb,verbatim,vmargin}
\usepackage[bookmarks=true]{hyperref}   

\theoremstyle{plain}
\newtheorem{thm}{Theorem}[section]
\newtheorem{lem}[thm]{Lemma}
\newtheorem{cor}[thm]{Corollary}
\newtheorem{prop}[thm]{Proposition}

\theoremstyle{definition}

\newtheorem{rem}[thm]{Remark}

\title[Remarks on maximal regularity]
{Remarks on maximal regularity}
\author{Pascal Auscher} \author{Andreas Axelsson} 
\address{Pascal Auscher, Universit\'e  Paris-Sud 11, UMR du CNRS 8628, 91405 Orsay Cedex, France and Centre for Mathematics and its Applications, Mathematical Sciences Institute, Australian National University, Canberra ACT 0200, Australia}
\email{pascal.auscher@math.u-psud.fr}
\address{Andreas Axelsson, Matematiska institutionen, Stockholms universitet, 106 91 Stockholm, Sweden}
\email{andax@math.su.se}

\mathchardef\semic="303B
\newcommand{\Mcc}{{M\raise.55ex\hbox{\lowercase{c}}}}

\newcommand{\R}{{\mathbf R}}

\newcommand{\mH}{{\mathcal H}}

\newcommand{\nul}{\textsf{N}}
\newcommand{\ran}{\textsf{R}}
\newcommand{\dom}{\textsf{D}}

\newcommand{\clos}[1]{\overline{#1}}

\newcommand{\barint}{\mbox{$ave \int$}}

\newcommand{\qe}[1]{\int_0^\infty\|#1\|^2\,\frac{dt}t}

\newcommand{\triple}[1]{|||#1|||}
\newcommand{\mM}{{\mathcal M}}

\def\barint_#1{\mathchoice
            {\mathop{\vrule width 6pt
height 3 pt depth -2.5pt
                    \kern -8.8pt
\intop}\nolimits_{#1}}%
            {\mathop{\vrule width 5pt height
3 pt depth -2.6pt
                    \kern -6.5pt
\intop}\nolimits_{#1}}%
            {\mathop{\vrule width 5pt height
3 pt depth -2.6pt
                    \kern -6pt
\intop}\nolimits_{#1}}%
            {\mathop{\vrule width 5pt height
3 pt depth -2.6pt
          \kern -6pt \intop}\nolimits_{#1}}}

\usepackage{color}

\definecolor{gr}{rgb}   {0.,   0.8,   0. } 
\definecolor{bl}{rgb}   {0.,   0.5,   1. } 
\definecolor{mg}{rgb}   {0.7,  0.,    0.7}

\begin{document}

\begin{abstract}
We prove weighted estimates for the maximal regularity operator. Such estimates were motivated by boundary value problems. We take this opportunity to study a class of weak solutions to the abstract Cauchy problem. We also give  a new proof of maximal regularity for closed and maximal accretive operators following from Kato's inequality for fractional powers and almost orthogonality arguments. 
\end{abstract}



\maketitle

\begin{center}
{\it In honour of H. Amann's birthday}
\end{center}

\section{Weighted estimates for the maximal regularity operator}

Assume $-A$ is a densely defined, closed linear operator,  generating a  bounded analytic semigroup $\{e^{-zA}, |\arg z\,| < \delta \}$, $0<\delta <\pi/2$,   on a Hilbert space $\mH$. 
Equivalently, $A$ is sectorial of type $\omega(A)=\pi/2- \delta $. 
Let $\dom(A)$ denote its domain. The maximal regularity operator is defined by the formula
$$\mM_{+}f(t)=\int_{0}^t Ae^{-(t-s)A} f(s)\, ds.$$
This operator is associated to the forward abstract evolution equation
$$ \dot u(t)+ Au(t)=f(t), t>0; \quad u(0)=0$$
as for appropriate $f$, $Au(t)=\mM_{+}f(t)$. An estimate on $\mM_{+}f$ in the same space as $f$ gives therefore bounds on $\dot u$ and $Au$ separately. See Section \ref{Cauchy}.

The integral defining $\mM_{+}f$ converges strongly in $\mH$ for  each  $f\in L^2(0,\infty; dt, \dom(A))$ and $t>0$. 
The estimate $\|Ae^{-(t-s)A}\|\le C(t-s)^{-1}$ following from the analyticity of the semigroup shows that the integral is singular if one only assumes $f(s) \in \mH$.  The maximal regularity operator is an example of a singular integral operator with  operator-valued kernel. The celebrated theorem by de Simon \cite{deS} asserts

\begin{thm}\label{thm1}  Assume $-A$ generates a bounded holomorphic semigroup in $\mH$. 
The operator $\mM_{+}$, initially defined on $L^2(0,\infty; dt, \dom(A))$, extends to a bounded operator
on $L^2(0,\infty; dt, \mH)$.
\end{thm}

Motivated by boundary value problems for some second order elliptic equations, we proved in 
\cite{AA} the following result.

\begin{thm} \label{AA}Assume $-A$ generates a bounded holomorphic semigroup in $\mH$ and  furthermore that $A$ has bounded holomorphic functional calculus, then  $\mM_{+}$, initially defined on $L^2_{c}(0,\infty; dt, \dom(A))$, extends to a bounded operator on $L^2(0,\infty; t^\beta dt, \mH)$ for all $\beta \in (-\infty, 1)$.
\end{thm}

The proof given there uses the operational calculus defined in the thesis of Albrecht \cite{A}.  It used as an assumption that $A$ has  bounded holomorphic functional calculus as defined by McIntosh \cite{Mc}. Under this assumption estimates of integral operators more general than the maximal regularity operator, with operator-kernels defined through
functional calculus of $A$, were proved and gave other useful informations  to understand also the case $\beta=1$ needed for the boundary value problems. However, not all generators of bounded analytic semigroups have a bounded holomorphic functional calculus. (See \cite{McY}, and  Kunstmann and Weis \cite[Section 11]{KW} for a list of equivalent conditions.) So if we only  consider the maximal regularity operator,  it is natural to ask whether one can drop the assumption on bounded holomorphic functional calculus in Theorem \ref{AA}. It is indeed the case and as we shall see the proof is extremely simple assuming we know Theorem \ref{thm1}.

\begin{thm}\label{thm1.3} Let $-A$ be the generator of a bounded analytic semigroup on $\mH$. Then $\mM_{+}$, initially defined on $L^2_{c}(0,\infty; dt, \dom(A))$, extends to a bounded operator on $L^2(0,\infty; t^\beta dt, \mH)$ for  all $\beta \in (-\infty, 1)$.
\end{thm}

The subscript $_{c}$ means with compact support in $(0,\infty)$.  Set $\triple{f(t)}^2=\qe{f(t)}$ (we leave in the $t$-variable in the notation for convenience).  As we often use it, we recall the following simplified version of Schur's lemma: if $U(t,s)$, $s,t>0$, are  bounded linear operators on $\mH$ with bounds $\|U(t,s)\| \le h(t/s) $ and 
$C=\int_{0}^\infty h(u) \frac{du} u <\infty$, then 
$$
\triple {\int_{0}^\infty U(t,s) f(s) \frac {ds}s } \le C \triple{f(s)}.
$$

\begin{proof}[Proof of Theorem 1.3]
 Let $\beta <1$. For $\beta=0$, this is Theorem \ref{thm1}. Assume $\beta\ne 0$ and set  $\alpha=\beta/2$. Observe that
$$ \|\mM_{+}f(t)\|_{L^2(t^\beta dt, \mH)}= \| { t^{\alpha}\mM_{+}f(t)}\|_{L^2(dt, \mH)}.$$
We have, with $f_{\alpha}(s)=s^{\alpha}f(s)$,  
$$ { t^{\alpha}\mM_{+}f(t)} = \mM_{+}(f_{\alpha})(t) + \int^{t}_{0} Ae^{-(t-s)A} (t^\alpha- s^\alpha) f(s)\, ds.
$$
For the first term apply Theorem \ref{thm1}. For the second, write
$$
\left\| \int^{t}_{0} Ae^{-(t-s)A} (t^\alpha- s^\alpha) f(s)\, ds
\right\|_{L^2(dt, \mH)}= \triple {\int_{0}^\infty U(t,s) g(s) \frac {ds}s }
$$
with $g(s) = s^{1/2+\alpha}f(s)$ and $U(t,s) = Ae^{-(t-s)A} (t^\alpha- s^\alpha) s^{1/2 - \alpha} t^{1/2}$ for $s< t$ and 0 otherwise. Since $\triple{g(t)}=\|f\|_{L^2(t^\beta dt, \mH)}$, it remains to estimate the norm of $U(t,s)$ on $\mH$.  We have
$$
\| U(t,s)\| \le C \frac{|t^\alpha- s^\alpha|}{|t-s|} s^{1/2 - \alpha} t^{1/2}, \quad s<t.
$$
It is easy to see that it is on the order of $(s/t)^{1/2 - \max(\alpha,0)}$ as $s< t$. 
We conclude by applying Schur's lemma.
\end{proof}

Let
$$\mM_{-}f(t)= \int_{t}^\infty Ae^{-(s-t)A} f(s)\, ds.$$
This operator is associated to the backward abstract evolution equation 
$$ \dot v(t)- Av(t)=f(t), t>0; \quad v(\infty)=0$$
as for appropriate $f$, $Av(t)=-\mM_{-}f(t)$.

\begin{cor} Assume that $-A$ generates a bounded analytic semigroup on $\mH$.  
Then $\mM_{-}$, initially defined on $L^2_{c}(0,\infty; dt, \dom(A))$, extends to a bounded operator on $L^2(0,\infty; t^\beta dt, \mH)$ for  all $\beta \in (-1,\infty)$.
\end{cor}

\begin{proof}
Observe that the adjoint of $\mM_{-}$ in $L^2(0,\infty; t^\beta dt, \mH)$  for the duality defined by $L^2(0,\infty; dt, \mH)$ is $\mM_{+}$ in $L^2(0,\infty; t^{-\beta} dt, \mH)$ associated to $A^*$ and apply  Theorem  \ref{thm1.3}.
\end{proof}

We  next  show that the range of $\beta$ is optimal in both results.

\begin{thm} For any non zero $-A$ generating a bounded analytic semigroup on $\mH$ and $\beta\ge 1$,   $\mM_{+}$ is not bounded on $L^2(0,\infty; t^\beta dt, \mH)$ and $ \mM_{-}$ is not bounded on $L^2(0,\infty; t^{-\beta} dt, \mH)$.

\end{thm}

\begin{proof}  It  suffices to  consider $\mM_{-}$. Since $A\ne 0$, 
  $\clos{\ran(A)}$, the closure of the range of $A$,  contains non zero elements. As $\clos{\ran(A)} \cap \dom(A)$ is dense in it,  pick $u\in \clos{\ran(A)} \cap \dom(A)$, $u\ne 0$, and  set $f(t)=u$ for $1\le t \le 2$ and 0 elsewhere. Then $ f \in L^2_{c}(0,\infty; dt, \dom(A))$ and $f\in  L^2(0,\infty; t^{-\beta} dt, \mH)$ with $\|f(t)\|_{L^2(0,\infty; t^{-\beta} dt, \mH)}= c_{\beta}\|u\| <\infty$. For $t<1$, one has 
$$
\mM_{-}f(t)= (e^{-(1-t)A} -e^{-(2-t)A}) u,
$$
 which converges to $ (e^{-A} -e^{-2A}) u $ in $\mH$ when $t\to 0$. We claim that $(e^{-A} -e^{-2A}) u\ne 0$ so 
$$\|\mM_{-}f(t)\|_{L^2(0,\infty; t^{-\beta} dt, \mH)}^2 \ge \int_{0}^1 \|(e^{-(1-t)A} -e^{-(2-t)A}) u\|^2 \, \frac {dt}{t^\beta} =\infty.$$

To  prove the claim, we argue as follows. Assume it is 0, then $e^{-2A}u= e^{-A}u$  so that an iteration yields $e^{-nA}u= e^{-A}u$ for all  integers $n\ge 2$. If $n \to \infty$, $e^{-nA}u$ tends to 0 in $\mH$
because $u \in \clos{\ran(A)}$. Thus $e^{-A}u=0$ and it follows that $e^{-tA}u=e^{-(t-1)A}e^{-A}u=0$ for all $t>1$. The
analytic function $z \to e^{-zA}u$ is thus identically 0 for $|\arg z\, |<\delta $. On letting $z \to 0$, we get
$u=0$ which is a contradiction.  
\end{proof}

We have seen that $\mM_{-}$ cannot map $L^2(0,\infty; t^{-1} dt, \mH)$ into itself and that it seems due to the behavior of $\mM_{-}f(t)$ at $t=0$ for some $f$. We shall make this precise and general: under a further assumption on $A$ which we introduce next, we define $\mM_{-} : L^2(0,\infty; t^{-1} dt, \mH) \to L^2_{loc}(0,\infty; dt, \mH)$ and show that   controlled behavior at 0 of $\mM_{-}f$ guarantees $\mM_{-}f \in L^2(0,\infty; t^{-1} dt, \mH)$.

We begin by writing whenever $f\in L^2_{c}(0,\infty; dt, \dom(A))$ and denoting $f_{-1/2}(s)=s^{-1/2}f(s)$,
\begin{align*}
\mM_{-}f(t) - e^{-tA} \int_{0}^\infty Ae^{-sA}f(s)\, ds 
&= t^{1/2}\mM_{-}(f_{-1/2})(t) 
\\
&
 +\int_{t}^{2t} Ae^{-(s-t)A}(s^{1/2}- t^{1/2})s^{1/2} f(s)\, \frac {ds}s 
 \\
 &
 +  \int^{\infty}_{2t} A(e^{-(s-t)A}-e^{-(s+t)A} )(s^{1/2}- t^{1/2})s^{1/2} f(s)\,  \frac {ds}s  
 \\
 &
 - \int^{\infty}_{2t} Ae^{-(s+t)A}   t^{1/2}s^{1/2} f(s)\,  \frac {ds}s   
 \\
 &-   \int_{0}^{2t} Ae^{-(s+t)A}s f(s) \,  \frac {ds}s. 
\end{align*}
 The right hand side is seen to belong to $L^2(0,\infty; t^{-1} dt, \mH)$ with an estimate 
 $C\triple{f(s)}$ using Theorem \ref{thm1} for the first term  and Schur's lemma for the other four terms. Hence, by density, the right hand side defines a bounded linear operator $\widetilde\mM_{-}$ on $L^2(0,\infty; t^{-1} dt, \mH)$. Also, the integral $\int_{0}^\infty Ae^{-sA}f(s)\, ds$ is defined as a Bochner  integral in $\mH$ whenever $f \in L^2_{c}(0,\infty; dt, \mH)$. Thus, by density of $\dom(A) $ in $\mH$, one can set for  $f \in L^2_{c}(0,\infty; dt, \mH)$,
\begin{equation}
\label{def}
 \mM_{-}f(t) : = \widetilde\mM_{-}f(t) + e^{-tA} \int_{0}^\infty Ae^{-sA}f(s)\, ds \quad \mathrm{in}\
 L^2_{loc}(0,\infty; dt, \mH).
\end{equation}
Let $E$ be the space
 of $f\in L^2(0,\infty; t^{-1} dt, \mH)$ such that the integrals $\int_{\delta }^R Ae^{-sA}f(s)\, ds$ converge weakly in $\mH$ as $\delta \to 0$ and $R\to \infty$.  Then  the  above equality extends to $f\in E$. Assuming, in addition, that $A^*$ satisfies the quadratic estimate \begin{equation}
\label{qe}
\triple{ sA^*e^{-sA^*}h} \le C\|h\|_{\mH} \quad \mathrm{for \ all}\ h\in \mH,
\end{equation} we have $E=L^2(0,\infty; t^{-1} dt, \mH)$. Indeed, 
for  all $f \in L^2(0,\infty; t^{-1} dt, \mH)$ and  $h\in \mH$,
\begin{equation}
\label{weak}
\int^{\infty}_{0}  \left|(sA e^{-sA}f(s), h) \right| \frac {ds}s  \le \triple{f(s)}\, \triple{ sA^*e^{-sA^*}h} \lesssim \triple{f(s)} \, \|h\|_{\mH}
\end{equation}
and the weak convergence of the truncated integrals follows easily.  Thus, the right hand side of \eqref{def} makes sense for all $f\in L^2(0,\infty; t^{-1} dt, \mH)$ under \eqref{qe} and this defines $\mM_{-}f$. Moreover, it follows from \eqref{weak} that 
\begin{equation}
\label{bound}
\sup_{\tau>0} \frac 1 \tau \int_{\tau}^{2\tau} \|\mM_{-}f(t)\|_{\mH}^2 \, dt \le C\triple{f(s)}^2.
\end{equation}
Then remark that 
\begin{equation}
\label{limit}
\lim_{\tau\to 0} \frac 1 \tau \int_{\tau}^{2\tau} \mM_{-}f(t) \, dt = \int^{\infty}_{0}  A e^{-sA}f(s) \, {ds} \quad \mathrm{in}\  \mH,
\end{equation}
as the corresponding limit for $\widetilde\mM_{-}f$ is 0 and  $e^{-tA} \to I$ strongly when $t\to 0$. 

All this yields  the following result.
\begin{prop} Let $-A$ be the generator of a bounded analytic semigroup in $\mH$ and assume that the quadratic estimate \eqref{qe} holds for $A^*$. Then \eqref{def} defines  $\mM_{-}f \in L^2_{loc}(0,\infty; dt, \mH)$ with estimates  \eqref{bound}  and limit \eqref{limit} for all $f\in L^2(0,\infty; t^{-1} dt, \mH)$.  In particular, 
$$\mM_{-}f\in L^2(0,\infty; t^{-1} dt, \mH)$$
 if  and only if 
$$
\lim_{\tau\to 0} \frac 1 \tau \int_{\tau}^{2\tau} \mM_{-}f(t) \, dt =0.
$$
The last condition defines a closed subspace of $L^2(0,\infty; t^{-1} dt, \mH)$ and there is a constant $C$ such that for all $f$ in this subspace
$$
\|\mM_{-}f(t)\|_{L^2(0,\infty; t^{-1} dt, \mH)}\le C \|f(t)\|_{L^2(0,\infty; t^{-1} dt, \mH)}.
$$

\end{prop}

Note that \eqref{qe} holds if $A$ has bounded holomorphic functional calculus by McIntosh's theorem \cite{Mc}.

\begin{rem}
For $\mM_{+}$, the analysis is not that satisfactory (for $\beta=1$). One can show similarly that
$$
\bigg\|\mM_{+}f(t) - Ae^{-tA}\int^{\infty}_{0}   e^{-sA}f(s) \, {ds}\bigg\|_{L^2(0,\infty; tdt, \mH)}\le C \|f(t)\|_{L^2(0,\infty; tdt, \mH)}
$$
provided $f \in L^2_{c}(0,\infty; dt, \dom(A))$. If the quadratic estimate \eqref{qe} holds for $A$, this allows to extend $\mM_{+}$ to the space $\{f \in L^2_{loc}(0,\infty; dt, \mH) ; \int^{\infty}_{0}   e^{-sA}f(s) \, {ds} \mathrm{\ converges\ weakly\ in \ } \mH\}$. However, there is no simple description of this space. 
\end{rem}

\section{Applications to the abstract Cauchy problem}\label{Cauchy}

In this section, we assume throughout that $-A$ generates a bounded analytic semigroup in $\mH$. 

 Let $f\in L^2_{loc}(0,\infty; dt, \mH) $. We say that $u$ is a weak solution  to $ \dot u(t)+ Au(t)=f(t), t>0,$   if $u \in L^2_{loc}(0,\infty; dt, \mH) $,  
\begin{equation}
\label{sup}
\sup_{0<\tau<1} \frac 1 \tau \int_{\tau}^{2\tau}  \left \|u(s) \right \|_{\mH} \, ds <\infty
\end{equation}
and for all $\phi \in C^1_{c}(0,\infty; \mH) \cap C^0_{c}(0,\infty; \dom(A^*))$,
\begin{equation}
\label{ODE}
\int_{0}^\infty (u(s), -\dot \phi(s)+A^*\phi(s)) \, ds = \int_{0}^\infty (f(s), \phi(s)) \, ds.
\end{equation}

The notion of weak solution here differs from the one in Amann's book \cite[Chapter 5]{Am}  called weak $L_{p,loc}$ solution ($p\in [1,\infty]$) specialized to $p=2$. We assume a uniform control through \eqref{sup} near $t=0$ and assume $\phi$ compactly supported in $(0,\infty)$ in \eqref{ODE}   instead of specifying the initial value at $t=0$ and taking  $\phi$ compactly supported in $[0,\infty)$ in \cite{Am}. 

\begin{lem} 
\label{v}Let $\beta\in (-\infty,1)$ and $f\in L^2(0,\infty; t^\beta dt, \mH)$. Then 
\begin{equation}
\label{X}
v(t)=\int_{0}^t
e^{-(t-s)A}f(s)\, ds
\end{equation} satisfies
\begin{enumerate}
  \item[(1)]  $v\in C^0([0,\infty); \mH)$ and  for all $t>0$,
  $\|v(t)\|^2 \le C t^{1-\beta} \int_{0}^t s^\beta\|f(s)\|^2 \, ds,
$ 
  \item[(2)]  $v$  is a weak solution to $ \dot u(t)+ Au(t)=f(t), t>0,$
  \item[(3)]  $Av(t)=\mM_{+}f(t)$ in $L^2_{loc}(0,\infty; dt, \mH) $, 
and 
$$\|\dot v(t)\|_{L^2(0,\infty; t^{\beta} dt, \mH)}+ \|Av(t)\|_{L^2(0,\infty; t^{\beta} dt, \mH)} \le C \|f(t)\|_{L^2(0,\infty; t^{\beta} dt, \mH)}.$$
Here, by $\mM_{+}$ we mean the bounded extension to $L^2(0,\infty; t^\beta dt, \mH)$.
\end{enumerate}

\end{lem}

\begin{proof} The inequality in (1) follows from the uniform boundedness of the semigroup and Cauchy-Schwarz inequality, and this shows that the integral defining $v(t)$ norm converges in $\mH$, thus infering  continuity on $[0,\infty)$,  and also  \eqref{sup}. To check \eqref{ODE}, it suffices to change order of integration and calculate. The equality $\mM_{+}f=Av$ is proved by duality against a $\phi$ as in \eqref{ODE} since such $\phi$ form a dense subspace in $L^2_{c}(0,\infty; dt, \mH) $. Finally, the inequalities in (3) are consequences of Theorem \ref{thm1.3}.
\end{proof}

We now state that all weak solutions have an explicit representation and a trace at $t=0$.  

\begin{prop}\label{representation} Let $\beta\in (-\infty,1)$ and $f\in L^2(0,\infty; t^\beta dt, \mH)$. Let $u$ be a weak solution to  $ \dot u(t)+ Au(t)=f(t), t>0.$ Then, there exists $h\in \mH$ such that 
\begin{equation}
\label{rep}
u(t)= e^{-tA}h + v(t) \quad \mathrm{in}\ L^2_{loc}(0,\infty; dt, \mH), 
\end{equation}
with $v$ defined by  \eqref{X}. In particular,   $t\mapsto u(t)$ can be redefined on a null set to be $C^0([0,\infty); \mH)$ with trace $h$ at $t=0$. \end{prop}

This immediately implies the following existence and uniqueness results.

\begin{cor}  Let $u_{0}\in \mH$. The initial value problem $ \dot u(t)+ Au(t)=0, t>0,$ with $\lim_{\tau\to 0} \frac 1 \tau \int_{\tau}^{2\tau} u(t)\, dt = u_{0}$ in $ \mH$,  has a unique  weak solution  given by $u(t)=e^{-tA}u_{0}$ for almost every $t>0$. In particular, up to redefining $t\mapsto u(t) $ on a null set, $u\in C^\infty(0,\infty; D(A))$ and is a strong solution. 
\end{cor}

\begin{cor}  Let $\beta\in (-\infty,1)$ and $f\in L^2(0,\infty; t^\beta dt, \mH)$. The initial value problem $ \dot u(t)+ Au(t)=f(t), t>0,$ with $ \lim_{\tau\to 0} \frac 1 \tau \int_{\tau}^{2\tau} u(t)\, dt = 0$ in $\mH$, has a unique weak solution given by $v$ defined by \eqref{X}, up to redefining $t\mapsto u(t)$ on a null set.

\end{cor}

\begin{proof}[Proof of Lemma \ref{representation}]  Define $\eta(s)$ to be the piecewise linear continuous function with support
$[1,\infty)$, which equals $1$ on $(2,\infty)$ and is linear on $(1,2)$.
Let $t>0$. For $0<\epsilon<t/4$ and $s>0$, let
$$
  \eta_\epsilon(t,s):= \eta(s/\epsilon)\eta( (t-s)/\epsilon).
$$
Let $\phi_0\in \mH$ be any boundary function, and choose 
$$\phi(s):= \eta_\epsilon(t,s) e^{-(t-s)A^*} \phi_0\in \mathrm{Lip}_{c}(0,\infty; \dom(A^*))$$  as test function (by approximating  $\eta_\epsilon(t,s) $ by a smooth function, this can be done). 
A calculation yields
\begin{align*}
-  \frac 1 \epsilon \int_{\epsilon}^{2\epsilon} \left(e^{-(t-s)A}u(s) ,  \phi_0\right)  \, ds +  \frac 1 \epsilon \int_{\epsilon}^{2\epsilon} & \left(e^{-sA}u(t-s),  \phi_0\right)\, ds
\\
&= 
   \int_0^\infty \left( \eta_\epsilon(t,s) e^{-(t-s)A} f(s), \phi_0\right)  \, ds
\end{align*}
 and since this is true for arbitrary $\phi_{0} \in \mH$ and $\eta_{\epsilon}$ has compact support, we deduce that
$$-  \frac 1 \epsilon \int_{\epsilon}^{2\epsilon} e^{-(t-s)A}u(s) \, ds  +  \frac 1 \epsilon \int_{\epsilon}^{2\epsilon} e^{-sA}u(t-s) \, ds= 
  \int_0^\infty\eta_\epsilon(t,s) e^{-(t-s)A} f(s)\, ds .
  $$
Now, we let $\epsilon \to 0$ as follows. First, $\eta_\epsilon(t,s) $ tends to the indicator function of $(0,t)$ so that the right hand side  is easily seen to converge to $v(t)$ in $\mH$ for any fixed $t>0$ by dominated convergence. Fix now $0<a<b<\infty$ and integrate in $t\in (a,b)$ the left hand side. Remark that $\frac 1 \epsilon \int_{a}^b\int_{\epsilon}^{2\epsilon} e^{-sA}u(t)\, dsdt$ converges to $\int_{a}^b u(t)\, dt $ in $\mH$. Substracting this quantity from the second term in the right hand side and using $u\in L^2_{loc}(0,\infty; \mH)$,  Lebesgue's theorem  yields
$$
\int_{a}^b\left \| \frac 1 \epsilon \int_{\epsilon}^{2\epsilon} e^{-sA}(u(t-s)-u(t))  \, ds\right\|_{\mH}^2 dt\le  \frac C \epsilon\int_{a}^b\int_{\epsilon}^{2\epsilon}  \left \| u(t-s)-u(t)\right\|_{\mH}^2   \, dsdt \to 0. 
$$
For the first term, using  $\|e^{-(t-s)A} - e^{-tA}\| \le Cs/t$ from analyticity and \eqref{sup}, one sees that 
\begin{equation}
\label{ana}
\left\|\frac 1 \epsilon \int_{\epsilon}^{2\epsilon} (e^{-(t-s)A}-e^{-tA})u(s) \, ds\right\|_{\mH} \to 0
\end{equation}
for each $t>0$. Thus 
$$
h_{\epsilon}(t):=  e^{-tA}h_{\epsilon}, \quad \mathrm{with}\quad h_{\epsilon}: = \frac 1 \epsilon \int_{\epsilon}^{2\epsilon} u(s) \, ds,
 $$
 has a  limit, say $h(t)$, in $L^2(a,b;\mH)$. The semigroup property yields 
 $h_{\epsilon}(t) = e^{-(t-\tau)A}  h_{\epsilon}(\tau)$ for all $t\ge\tau $. Thus, 
 $$
  \| h_{\epsilon}(t)-  h_{\epsilon'}(t)\|_{\mH} 
  \le \frac 1{b-a}\int_a^b \| e^{-(t-\tau)A} ( h_{\epsilon}(\tau)-  h_{\epsilon'}(\tau)) \|_\mH d\tau
  \le C \left( \int_a^b \| h_{\epsilon}(\tau)-  h_{\epsilon'}(\tau)\|^2_{\mH} d\tau\right)^{1/2},
$$
when $t>b$.
Hence, since $(a,b)$ is arbitrary, $h_{\epsilon}(t)$ converges in $\mH$ to $h(t)$ for each $t>0$. Thus, for any $\phi_{0} \in \mH$ and $t>0$, we have
$$
(h_{\epsilon}, e^{-tA^*}\phi_{0})= (h_{\epsilon}(t), \phi_{0}) \to (h(t),\phi_{0}).
$$
Since $(h_{\epsilon})_{\epsilon<1}$ is a bounded sequence in $\mH$ by \eqref{sup} and the elements $e^{-tA^*}\phi_{0}$, $ t>0$, $\phi_{0}\in \mH$, form a dense set of $\mH$, we infer that $h_{\epsilon}$ has a weak limit in $\mH$. Calling $h$ this weak limit we have $
(h, e^{-tA^*}\phi_{0})= (h(t),\phi_{0})$, hence
$h(t)=e^{-tA}h$ as desired. Summarizing, we have obtained $-e^{-tA}h+u(t)=v(t)$ in $L^2(a,b;\mH)$ for all $0<a<b<\infty$. 

Thus, $u$ agrees almost everywhere with the continuous function  $t\mapsto v(t)+e^{-tA}h$  which has limit $h$ at $t=0$. 
\end{proof}

\begin{rem} The only time analyticity is used in this proof is in  \eqref{ana}. If we had incorporated the existence of an initial value as in \cite{Am} in our definition of a weak solution then  analogous proposition and corollaries would hold for all generators of  bounded $C^0$-semigroups.

\end{rem}

\section{A proof of maximal regularity via Kato's inequality for fractional powers}

There are many proofs of the de Simon's theorem, via Fourier transform or operational calculus, and various extensions to Banach spaces. We refer to   \cite[Section 1]{KW}.

Here, we wish to provide a proof using ``almost orthogonality arguments''  (Cotlar's lemma), and Kato's inequality for fractional powers \cite[Theorem 1.1]{K} which we recall for the reader's convenience.

\begin{thm}
 Let $A$ be  closed and maximal accretive. For any $0\le \alpha < 1/2$, the operators $A^\alpha$ and $A^{*\alpha}$ have same domains and satisfy
 \begin{equation}\label{comparability}
\| A^{*\alpha} f \| \le \tan \frac{\pi(1+2\alpha)}{4} \| A^\alpha f \|.
\end{equation}
If, moreover,  $A$ is injective then $A^\alpha A^{*-\alpha}$ extends to  a bounded operator on $\mH$ for $-1/2<\alpha<1/2$. 
\end{thm}

Maximal accretive means that ${\textrm Re} (Au , u) \ge 0$ for every $u \in \dom(A)$ and $(\lambda -A)^{-1} $ is bounded 
 whenever ${\text Re}\lambda < 0$. Note that \eqref{comparability} holds true with different constants for operators which are similar to a closed and maximal accretive  operator. Assume $A$ is sectorial of type $w(A)<\pi/2$ and injective. Le Merdy showed in \cite{LeM1} that  $A$ is similar to a maximal accretive operator if and only if $A$ has bounded imaginary powers (i.e. $A^{it}$ is bounded for all $t\in \R$). (See also \cite{LeM2} for a more general result and \cite{Sim} for explicit examples.) But, following earlier works of Yagi \cite{Y},  McIntosh showed in his seminal paper \cite{Mc} that  $A$ has bounded imaginary powers if and only if $A$ has a bounded holomorphic functional calculus.
(See \cite[Section 11]{KW}  for extensive discussions with historical notes.) So proving maximal regularity (i.e. Theorem \ref{thm1}) assuming maximal accretivity is the same as proving maximal regularity assuming bounded holomorphic functional calculus. Nevertheless, this direct argument below could be of interest.

\begin{proof}[Proof of Theorem \ref{thm1} under further assumption of maximal accretivity]
Since $Ae^{-(t-s)A}$ \break anni\-hilates $\nul(A)$, the null space of $A$, we may assume $g(s)\in \clos{\ran(A)}$ for all $s>0$. Alternately, we may factor out the null space of $A$ and assume that $A$ is injective, which we do ($A$ is sectorial, so $\mH$ splits topologically as $\nul(A) \oplus \clos{\ran(A)}$).

Then one can write 
$g(s)=\int_{0}^\infty uAe^{-uA}  g(s) \frac{du}u$ and so we have the representation of $\mM_{+}$ as
$$\mM_{+}g(t) = \int_{0}^\infty  ( T_{u}g)(t)\frac{du}u, \quad {\rm  with} \  (T_{u}g)(t)= \mM_{+}(uAe^{-uA}g)(t).
$$
By Cotlar's lemma (see \cite[Chapter VII]{Ste}) it is enough to show in operator norm 
on $L^2(0,\infty; \mH)$ that 
$\|T_{u}T_{v}^*\| +\|T_{u}^*T_{v}\|  \le h(u/v)$  with  $C=\int_{0}^\infty h(x) \frac{dx}x <\infty$
to conclude that $\mM_{+}$ is bounded on $L^2(0,\infty; \mH)$ with norm less than or equal to $C$.
We show that for all $\alpha\in (0,1/2)$ one can take $
h(x)= C_{\alpha} \min \left( x^\alpha, x^{-\alpha}\right).
$

We begin with $T_{u}T_{v}^*$ for fixed $(u,v)$. Since $\|T_{u}T_{v}^*\| = \|T_{v}T_{u}^*\|$, we may assume $u\le v$. A computation yields
$$
(T_{u}T_{v}^*)(g)(t)= \int_{0}^\infty K_{(u,v)}(t,\tau)g(\tau)\, d\tau
$$
where 
$$
K_{(u,v)}(t,\tau) = \int_{0}^{\min(t,\tau)} uA^2e^{-(t-s+u)A}   v{A^*}^2e^{-(\tau-s+v)A^*} \, ds.
$$ 
We turn to  estimate the operator norm on $\mH$ of $K_{(u,v)}(t,\tau)$ for fixed $(t, \tau)$. (Recall we fixed $(u,v)$ with $u\le v$.)   Since $A$ is maximal accretive and injective, we have  $\|A^\alpha {A^*}^{-\alpha}\| \le C(\alpha)$ for $\alpha\in (0,1/2)$.
So we write
$$
uA^2e^{-(t-s+u)A}   vA^{*2}e^{-(\tau-s+v)A^*} = uA^{2-\alpha}  e^{-(t-s+u)A} (A^\alpha {A^*}^{-\alpha}) v{A^*}^{(2+\alpha)}e^{-(\tau-s+v)A^*},
$$
and by analyticity the operator norm on $\mH$ is bounded by constant times $a(s)b(s)$ with
$$
a(s)=  \frac{ u}{(t-s+u)^{2-\alpha}} , \quad b(s) = \frac{ v}{(\tau-s+v)^{2+\alpha}}.
 $$
 Plug this estimate into the integral. If $t\le \tau$, bound $b(s)$ by $b(t)$ and get  
 $$
 \|K_{(u,v)}(t,\tau) \| \le C  u^{\alpha} b(t) = C (u/v)^\alpha \frac{ v^{1+\alpha}}{(\tau-t+v)^{2+\alpha}}.$$
 If $\tau\le t$,  bound $a(s)$ by $a(\tau)$  and get
 $$
  \|K_{(u,v)}(t,\tau) \| \le C a(\tau)  v^{-\alpha}= C (u/v)^\alpha \frac{ u^{1-\alpha}}{(t-\tau+u)^{2-\alpha}}. 
$$
It follows that 
$$
\sup_{\tau>0}\int_{0}^\infty  (\|K_{(u,v)}(t,\tau) \| +  \|K_{(u,v)}(\tau,t) \|)\, dt \le C(u/v)^\alpha.
$$
By Schur's lemma we obtain  $\|T_{u}T_{v}^*\| \le C(u/v)^\alpha$ when $u\le v$. 

We now turn to  estimate $T_{u}^*T_{v}$. By symmetry under taking adjoints again, it is enough to assume $u\le v$. We obtain
$$
(T_{u}^*T_{v})(g)(t)= \int_{0}^\infty \tilde K_{(u,v)}(t,\tau)g(\tau)\, d\tau
$$
where 
$$
\tilde K_{(u,v)}(t,\tau) = \int^{\infty}_{\max(t,\tau)} u{A^*}^2e^{-(s-t+u)A^*}   v{A}^2e^{-(s-\tau+v)A} \, ds.
$$ 
 This time we use the bound $\|{A^*}^\alpha {A}^{-\alpha}\| \le C(\alpha)$ for $\alpha\in (0,1/2)$
to obtain, if $\tau \le t$,
$$
 \|\tilde K_{(u,v)}(t,\tau) \| \le  C (u/v)^\alpha \frac{ v^{1+\alpha}}{(\tau-t+v)^{2+\alpha}}$$
 and if  $t \le \tau$,  
 $$
  \|\tilde K_{(u,v)}(t,\tau) \| \le  C (u/v)^\alpha \frac{ u^{1-\alpha}}{(t-\tau+u)^{2-\alpha}} .
$$
So, 
$$
\sup_{\tau>0}\int_{0}^\infty  (\|\tilde K_{(u,v)}(t,\tau) \| +  \|\tilde K_{(u,v)}(\tau,t) \|)\, dt \le C(u/v)^\alpha
$$
and by Schur's lemma, $\|T_{u}T_{v}^*\| \le C(u/v)^\alpha$ when $u\le v$.
\end{proof}

As Kato's inequality holds for all $\alpha \in (-1/2,1/2)$, the argument above  can be used to prove that $\mM_{+}$ is bounded on $L^2(0,\infty; t^\beta dt, \mH)$ but for $\beta \in (-1,1)$. We leave details to the reader.  

\

We thank Alan McIntosh for discussions on  the topic of this short note.

\end{document}